\pdfoutput=1
\RequirePackage{ifpdf}
\ifpdf % We are running pdfTeX in pdf mode
\documentclass[pdftex]{sigma}
\else
\documentclass{sigma}
\fi

\numberwithin{equation}{section}

\numberwithin{theorem}{section}
\numberwithin{proposition}{section}
\numberwithin{lemma}{section}
\numberwithin{corollary}{section}
\numberwithin{definition}{section}
\numberwithin{example}{section}
\numberwithin{remark}{section}
\numberwithin{note}{section}

\newcommand{\RR}{\ensuremath{\mathbb{R}}}

\DeclareMathOperator\Div{div}

\begin{document}

\allowdisplaybreaks

\renewcommand{\PaperNumber}{035}

\FirstPageHeading

\ShortArticleName{First Integrals of Vector Fields}

\ArticleName{A Note on the First Integrals of Vector Fields\\ with Integrating Factors and Normalizers}

\Author{Jaume LLIBRE~$^\dag$ and Daniel PERALTA-SALAS~$^\ddag$}

\AuthorNameForHeading{J.~Llibre and D.~Peralta-Salas}

\Address{$^\dag$~Departament de Matem\`{a}tiques,
Universitat Aut\`{o}noma de Barcelona,\\
\hphantom{$^\dag$}~08193 Bellaterra, Barcelona,
Catalonia, Spain}
\EmailD{\href{mailto:jllibre@mat.uab.cat}{jllibre@mat.uab.cat}}
\URLaddressD{\url{http://www.gsd.uab.es/personal/jllibre}}

\Address{$^\ddag$~Instituto de Ciencias Matem\'{a}ticas, Consejo Superior de Investigaciones Cient\'ificas,\\
\hphantom{$^\ddag$}~C/ Nicol\'{a}s Cabrera 13-15, 28049 Madrid, Spain}
\EmailD{\href{mailto:dperalta@icmat.es}{dperalta@icmat.es}}
\URLaddressD{\url{http://www.icmat.es/dperalta}}

\ArticleDates{Received February 16, 2012, in f\/inal form June 12, 2012; Published online June 14, 2012}

\Abstract{We prove a suf\/f\/icient condition for the existence of explicit f\/irst
integrals for vector f\/ields which admit an integrating factor. This
theorem recovers and extends previous results in the literature on
the integrability of vector f\/ields which are volume preserving and
possess nontrivial normalizers. Our approach is geometric and
coordinate-free and hence it works on any smooth orientable
manifold.}


\Classification{34C05; 34A34; 34C14}

\section{Introduction and statement of the main result}\label{s1}

Let $M$ be a smooth orientable manifold of dimension $n$ endowed
with a volume form $\Omega$. We denote by $C^\infty(M)$ and
$\mathfrak X^\infty(M)$ the spaces of smooth real-valued functions
and vector f\/ields on~$M$ respectively. We shall only work in the
smooth ($C^\infty$) category, although most of our results can be
easily adapted to lower regularity.

This letter is focused on the integrability properties of vector
f\/ields. More precisely, we are interested in proving some suf\/f\/icient
conditions which imply the existence of explicit f\/irst integrals
(i.e.\ computable). We recall that a function $H\in C^\infty(U)$,
where $U\subseteq M$ is an open subset of~$M$, is a {\em
first integral} of the vector f\/ield~$X$ in~$U$ if $L_X(H)=0$, which implies
that $H$ is constant along the integral curves of $X$. As usual,
$L_X$ denotes the Lie derivative with respect to the vector f\/ield~$X$. We recall that an equivalent way of writing the f\/irst integral
condition is $L_X(H)=i_X(dH)=X(H)=0$, where $d$ is the exterior
derivative and $i$ denotes de contraction operator. In this paper we are interested in global or semi-global f\/irst integrals in the sense that the open set $U\subseteq M$ where the f\/irst integral is well def\/ined is known a priori.

The importance of the existence of a non-constant f\/irst integral lies in the fact
that the trajectories of the vector f\/ield~$X$ leave invariant the
level sets of the function~$H$, and hence this is a strong
constraint on the dynamical behavior of the vector f\/ield (e.g.\ it
prevents from the existence of ergodic trajectories on~$M$). Unfortunately, it is
generally very dif\/f\/icult to compute a f\/irst integral (whose
existence is, in fact, a rather non-generic phenomenon). This is the
reason why it is necessary to introduce some auxiliary objects to
study the integrability properties of a vector f\/ield.

There are two important tools which can be used to analyze the
existence of f\/irst integrals: integrating factors and normalizers.
We say that a function $f\in C^\infty(M)$ is an {\em integrating
factor} of the vector f\/ield $X$ (with respect to a volume form
$\Omega$) if $\Div_\Omega(fX)=0$. Recall that the divergence of a
vector f\/ield $X$ with respect to a volume form~$\Omega$,
$\Div_\Omega(X)$, is def\/ined by
$L_X(\Omega)=\Div_\Omega(X) \Omega$. A
vector f\/ield $Y\in \mathfrak X^\infty(M)$ is a {\em normalizer} of
$X$ if $[X,Y]=\lambda X$ for some function $\lambda\in C^\infty(M)$ (the notion of normalizer was introduced by Lie and Engel, see e.g.\ the recent translation~\cite[Theorem~20, Chapter~8]{Me10}).
Remember that $[\cdot,\cdot]$ stands for the Lie bracket, which is def\/ined as
$[X,Y]=L_X(Y)$. Normalizers are also called orbital symmetries, and
in the case that $[X,Y]=0$, it is usually said that the vector f\/ield
$Y$ is a centralizer (or symmetry) of~$X$.

There is an abundant literature on the connection between f\/irst
integrals, integrating factors and normalizers. The fact that the
existence of normalizers gives rise to computable f\/irst integrals in terms of quadratures
is well known since Lie's works, see e.g.\ the monographs on the
subject by Olver~\cite{Ol93} and Bluman--Anco~\cite{BA02} or the comprehensive geometric account of reduction and integrability by Sherring and Prince~\cite{SP92}. Some recent results on the relationship between symmetries and f\/irst integrals can be consulted in~\cite{CGP01,MN03,PS08,Pr09} and references therein. On the other hand, the
existence of integrating factors is not generally enough to ensure
the existence of f\/irst integrals, except for the case that $M$ has
dimension $2$~\cite[Chapter~2.11]{Gor}. In the case that $M=\RR^n$ and the vector f\/ield $X$ is polynomial, some algebraic algorithms have been proposed to compute polynomial integrating factors and elementary f\/irst integrals~\cite{PS83,W00}. For volume-preserving vector
f\/ields, i.e.\ $\Div_\Omega(X)=0$ (a particular case of vector f\/ields admitting an integrating
factor) which have volume-preserving normalizers, there are several
procedures to compute explicit f\/irst integrals using quadratures~\cite{HM98,H04, MW94}.
We also note that integrating factors are also directly related to
normalizers, and in fact can be obtained from them under suitable
hypotheses, see e.g.~\cite{GG10}.

In this letter we follow a totally dif\/ferent approach, which combines the
existence of integrating factors with a suitable extension of the
notion of normalizer (see $(i)$ below). The main dif\/ference with the approaches explained in the previous paragraph is that we provide an explicit formula for the f\/irst integral of the vector f\/ield in terms of the normalizer and the integrating factor, and there is no need of quadratures, reduction of order or algebraic manipulations to get this f\/irst integral. Previous results in this direction were obtained in~\cite{G77,GP00,Ho92} for volume preserving vector f\/ields with centralizers (see Section~\ref{S.ex} for more details). Let us now state the main theorem of this letter.

\begin{theorem}\label{t1}
Let $X$ be a smooth vector field on a manifold $M$. Assume that the
following three conditions hold:
\begin{enumerate}\itemsep=0pt
\item[$(i)$] There exists a vector field $Y\in \mathfrak X^\infty(M)$
and two functions $\lambda, \mu\in C^\infty(M)$ such that $[X,Y]=
\lambda X+\mu Y$.

\item[$(ii)$] $X$ admits an integrating factor $f\in C^\infty(M)$.

\item[$(iii)$] The function $f\mu$ is an integrating factor of $Y$.
\end{enumerate}
Then the function
\begin{gather}\label{e0}
H:=\frac{L_Y f}{f}+ \Div_\Omega(Y)-\lambda
\end{gather}
is a smooth first integral of the vector field $X$ in the open set
$U:=M\backslash Z_f$, where $Z_f:=\{p\in M: f(p)=0\}$.
\end{theorem}

\begin{remark}
The theorem holds true if $X$, $Y$ and $f$ are of class $C^2$; in this case the f\/irst integral~$H$ is of class $C^1$. On the other hand, if the zero set $Z_f$ has empty interior the f\/irst integral $H$ is def\/ined in an open and dense subset of $M$. This is the case, e.g.\ if the integrating factor $f$ is analytic ($C^\omega$).
\end{remark}

Let us brief\/ly explain the geometric meaning of the f\/irst integral derived in Theorem~\ref{t1}, cf.\ equation~\eqref{e0}. We are very grateful to one of the referees for pointing out to us the argument that we sketch below. It is obvious that the function $g:=Hf$ is a {\em second integrating factor} of the vector f\/ield $X$; conversely, if $X$ admits two integrating factors $f$ and $g$ it is ready to check that the function $I:=g/f$ is a f\/irst integral of $X$, but it does not need to be the same as the f\/irst integral $H$ computed in Theorem~\ref{t1} (see e.g.\ Example~\ref{exvp} in Section~\ref{s2}). The point is that under the assumptions of Theorem~\ref{t1} it is possible to compute a second integrating factor $g$ using the vector f\/ield $Y$ and the integrating factor~$f$, such that the function $g/f$ is precisely the f\/irst integral $H$ given by equation~\eqref{e0}. For example, in the particular case that~$Y$ is a symmetry of $X$ (i.e.\ $[X,Y]=0$), then $g:=L_Y f$ is a second integrating factor of $X$, thus implying that $H:=\frac{L_Y f}{f}$ is a f\/irst integral of the vector f\/ield $X$. In the general case, one can look for a second integrating factor of the form $g:=L_Y f+\hat g$, and Theorem~\ref{t1} provides some suf\/f\/icient conditions to compute the unknown function $\hat g$. This is the main idea behind the proof of Theorem~\ref{t1}, which is given in Section~\ref{s2}.

Let us also observe that the f\/irst integral $H$ def\/ined in \eqref{e0} could be trivial (a
constant), in which case our theorem does not provide any new
information about the vector f\/ield $X$. The fact that this case can
indeed happen (even when $X$ has non-trivial f\/irst integrals) is illustrated in
Section~\ref{s2}, thus showing that our way of producing f\/irst
integrals is by no means the most general one. Finally, in
Section~\ref{S.ex} we give some applications and examples of our
main result and show how it generalizes some previous theorems in
the literature.

\section{Proof of Theorem~\ref{t1} and an example}\label{s2}

First we shall prove Theorem~\ref{t1}.

\begin{proof}[Proof of Theorem~\ref{t1}]
Our method of proof, which is coordinate-free, makes use of standard
tools from dif\/ferential geometry, the reader who is not familiar
with this dif\/ferential geometric setting can consult e.g.\
\cite[Chapter 2]{AM}. We f\/irst recall the following identity:
\begin{gather}\label{e1}
L_{[X,Y]}= L_XL_Y-L_YL_X .
\end{gather}

Using the properties of the Lie derivative and assumptions $(ii)$ and
$(iii)$ in Theorem~\ref{t1}, it is easy to check that
\begin{gather}\label{mu}
 L_{\lambda X}(f\Omega)=L_X(\lambda f\Omega)=L_X(\lambda)f\Omega ,\\
 L_{\mu Y}(f\Omega)=L_{f\mu Y}(\Omega)=\Div_\Omega(f\mu Y)
\Omega=0 .\label{mu2}
\end{gather}

We now compute $L_{[X,Y]}(f\Omega)$. On account of \eqref{mu} and
\eqref{mu2} and assumption $(i)$, the linearity of the Lie derivative
implies that
\begin{gather}\label{com1}
L_{[X,Y]}(f\Omega)=L_{\lambda X}(f\Omega)+L_{\mu Y}(f\Omega)=
L_X(\lambda)f\Omega .
\end{gather}

On the other hand, the identity~\eqref{e1}, assumption $(ii)$ and the fact that the
Lie derivative obeys the Leibniz rule allow to compute
$L_{[X,Y]}(f\Omega)$ in a dif\/ferent way, that is:
\begin{gather*}
 L_{[X,Y]}(f\Omega)=L_XL_Y(f\Omega)-L_YL_X(f\Omega)=L_X\left[\left(
\frac{L_Y(f)}{f}+\Div_\Omega(Y)\right)f\Omega\right]\\
 \hphantom{L_{[X,Y]}(f\Omega)}{} =\left[L_X\left(\frac{L_Y(f)}{f}+\Div_\Omega(Y)\right)\right]f\Omega .
\end{gather*}
Combining this equation with \eqref{com1}, we f\/inally conclude that
the function
\[
H:=\frac{L_Y(f)}{f}+\Div_\Omega(Y)-\lambda
\]
satisf\/ies the equation
\[
L_X(H)=0 ,
\]
which implies that $H$ is a f\/irst integral of the vector f\/ield $X$. Observe that $H$ is well def\/ined (and $C^\infty$) in the set $U:=M\backslash Z_f$, where $Z_f:=\{p\in M: f(p)=0\}$. This completes the proof of the theorem.
\end{proof}

\begin{remark}
Arguing exactly as in the proof of Theorem~\ref{t1} it is obtained that the function
\[
g:=L_Y f+f\Div_\Omega(Y)-\lambda f
\]
is a second integrating factor of the vector f\/ield $X$. This implies that $H:=g/f$ is a f\/irst integral of~$X$, thus providing an alternative approach to the proof of Theorem~\ref{t1} in the spirit of the explanation given at the end of Section~\ref{s1}: the f\/irst integral~$H$ follows from the existence of an appropriately chosen second integrating factor.
\end{remark}

We f\/inish this section by observing that, in the case that the
function $H$ is identically constant, it provides a trivial f\/irst
integral of $X$. The following example shows that this situation can
indeed happen, although it is not the general case as we shall see
in Section~\ref{S.ex}.

\begin{example}\label{exvp}
Let $X$ be a volume-preserving vector f\/ield having a
volume-preserving centra\-li\-zer~$Y$ which is (almost everywhere)
transverse to~$X$, that is
\[
[X,Y]=0 ,\qquad \Div_\Omega(X)=\Div_\Omega(Y)=0 ,\qquad
 \text{rank}(X,Y)=2\quad \text{a.e.\ in} \  M .
 \] It is
straightforward to check that assumptions $(i)$--$(iii)$ in
Theorem~\ref{t1} are satisf\/ied by the vector f\/ield~$X$. On the other
hand, the function~$H$ is identically zero, thus providing a trivial
f\/irst integral of~$X$. When the manifold~$M$ is $3$-dimensional, it
can be shown~\cite{HM98, MW94} that
\[
d(i_Xi_Y\Omega)=0 ,
\]
thus concluding, if the f\/irst De Rham cohomology group of $M$ is trivial, that is $H^1(M;\RR)=0$, that there exists a function
$I\in C^\infty(M)$ such that $dI=i_Xi_Y\Omega$. This function $I$ is
non-trivial because the vector f\/ield~$Y$ is transverse to~$X$, and
since $L_X(I)=i_XdI=i_Xi_Xi_Y\Omega=0$, it is a f\/irst integral of~$X$. Accordingly, the vector f\/ield~$X$ has a global f\/irst integral which is
not produced by the mechanism involved in Theorem~\ref{t1}.
\end{example}

\begin{remark}
Note that we have not assumed that the vector f\/ield $Y$ be (almost
everywhere) transverse to $X$ in the statement of Theorem~\ref{t1}.
If $Y$ is proportional to $X$, i.e.\ $Y=hX$ for some function $h$,
then we have that $\lambda=L_X(h)$, and being $\Div_\Omega(X)=
\frac{-L_X(f)}{f}$, it easily follows that $H=0$. We conclude that,
in order to obtain a non-trivial f\/irst integral $H$ with
Theorem~\ref{t1}, the vector f\/ield $Y$ must be transverse to $X$.
\end{remark}

\section{Applications and examples}\label{S.ex}

A particularly useful application of Theorem~\ref{t1} is when the
vector f\/ield~$Y$ is a normalizer of~$X$ (a~particular case of
assumption $(i)$ in the theorem) and $X$ admits an integrating factor.
In this case, a well known property of the normalizers of a vector f\/ield allows to obtain many
f\/irst integrals (not necessarily independent). The following
corollary is a generalization of a~theorem in~\cite{G77,GP00} which
states that $\Div_\Omega(Y)$ is a f\/irst integral of $X$ provided
that $[X,Y]=0$ and \mbox{$\Div_\Omega(X)=0$}.

\begin{corollary}\label{c0}
Assume that a vector field $X\in \mathfrak X^\infty(M)$ has a
normalizer $Y\in \mathfrak X^\infty(M)$, i.e.\ $[X,Y]=\lambda X$ for
some function $\lambda \in C^\infty(M)$. Then, if $X$ admits an
integrating factor~$f$, the function $H$ defined in~\eqref{e0} is a
first integral of~$X$. Moreover, the function
$H_k:=L_Y(L_Y(\cdots(L_Y(H))))$ is also a first integral of~$X$,
where the Lie derivative $L_Y$ is applied $k\geq 1$ times on~$H$.
\end{corollary}

\begin{proof}
It is obvious that the hypotheses of Theorem~\ref{t1} are fulf\/illed
with $\mu=0$, and therefore the function $H$ is a f\/irst integral of
$X$. If we prove that the function $H_1:=L_Y(H)$ is a f\/irst integral
of $X$ it easily follows by induction that any function $H_k$ is
also a f\/irst integral of $X$, so let us prove that this is the case
for $H_1$. Indeed,
\[
L_X(H_1)=L_X(L_Y(H))=L_{[X,Y]}(H)+L_Y(L_X(H))=\lambda L_X(H)=0 ,
\]
where we have used the identity~\eqref{e1} and that $L_X(H)=0$, thus proving the claim.
\end{proof}

Notice that at most $n-1$ independent f\/irst integrals can be
obtained with this procedure, so most of the functions $H_k$ will
not add any new information. On the other hand, it is relevant to
remark that the existence of an integrating factor for a vector
f\/ield is not enough in general to guarantee the existence of a f\/irst
integral, so it is necessary to assume further hypotheses, as in the
statement of Corollary~\ref{c0}.

In the following corollary we give a suf\/f\/icient condition in order
that a vector f\/ield be completely integrable, i.e.\ to admit $n-1$
independent f\/irst integrals (we say that two functions are
independent if they are functionally independent in an open and
dense subset of~$M$). Completely integrable vector f\/ields were
studied in~\cite{PS08} concerning their connection with the period
function and the existence of normalizers. Let us observe that this result is quite dif\/ferent from a theorem recently proved by Kozlov~\cite{Ko05}: on the one hand, we do not need to assume that the set of vector f\/ields $\{X,Y_1,\dots,Y_{n-2}\}$ generates a (solvable) Lie algebra and on the other hand the f\/irst integrals that we obtain are explicit (i.e.\ $\Div_\Omega(Y_i)$) and there is no need of using quadratures to compute them.

\begin{corollary}\label{c2}
Let $X\in\mathfrak X^\infty(M)$ be a vector field whose flow
preserves the volume form $\Omega$ and admits $n-2$ normalizers
$Y_1,\dots,Y_{n-2}$ so that $[X,Y_i]=c_iX$ with $c_i\in\mathbb R$,
$i=1,\dots,n-2$. Then, if the functions $\Div_\Omega(Y_i)$ are
independent, the vector field $X$ is completely integrable.
\end{corollary}

\begin{proof}
By assumption, $\Div_\Omega(X)=0$, and hence Corollary~\ref{c0} can
be applied with $f=1$ and $\lambda=c_i$ to imply the existence of
$n-2$ f\/irst integrals $H_i=\Div_\Omega(Y_i)$. Under the assumption that the functions $H_i$
are independent, it is standard that the pull-back of
the vector f\/ield $X$, which is volume-preserving, to the regular
level sets (of dimension $2$) of the map $(H_1,\dots,H_{n-2}):M\to
\RR^{n-2}$ admits an integrating factor with respect to the induced
volume form, and hence it is integrable~\cite[Chapter 2.11]{Gor}.
This def\/ines a f\/irst integral $H_{n-1}$ in any contractible domain of the regular level sets, thus proving the claim.
\end{proof}

\begin{remark}
According to Corollary~\ref{c2}, the vector f\/ield $X$ has $n-2$ f\/irst integrals given by $H_i=\Div_\Omega(Y_i)$ ($i=1,\dots,n-2$), and hence we can apply~\cite[Chapter 4.4, Theorem 13]{AKN97} to conclude that on each regular level set of the f\/irst integrals the f\/ield can be integrated by quadratures. Moreover, assuming that the regular level set is connected and compact and $X$ is non-vanishing, it follows that it is dif\/feomorphic to a torus and $X$ is orbitally conjugated to a linear f\/ield of (rational or irrational) frequency on the level set. Observe that, in order to use~\cite[Chapter 4.4, Theorem 13]{AKN97}, it is necessary to have $n-2$ f\/irst integrals, and Corollary~\ref{c2} precisely gives a suf\/f\/icient condition for this.
\end{remark}

Theorem~\ref{t1} can be applied to compute f\/irst integrals of Hamiltonian vector f\/ields with symmetries (note that a Hamiltonian vector f\/ield is divergence-free with respect to some volume form $\Omega$). In the following example, which is inspired by a paper of Hojman~\cite{Ho92}, we apply it to the motion of a particle in a spherically symmetric homogeneous potential to obtain the well known conservation laws of the energy and the angular momentum. We f\/ind this result rather surprising, as the method of proof is completely dif\/ferent from the classical one, thus showing the existence of a mechanism dif\/ferent from Noether's theorem to produce f\/irst integrals in Hamiltonian mechanics.

\begin{example}
Let $V(u)$ be a real valued function which is homogeneous of degree $\alpha\in \RR\backslash\{0\}$, that is $V(u):=u^\alpha$. Endowing $\RR^{2n}$ with Cartesian coordinates $(x_1,\dots,x_n,p_1,\dots,p_n)$, the Hamiltonian vector f\/ield which describes the dynamics of a particle under the action of the potential $V(x_1^2+\dots+x_n^2)$ is given by
\[
X=p_1\partial_{x_1}+\dots+p_n\partial_{x_n}-2x_1V'\partial_{p_1}-\dots-2x_nV'\partial_{p_n} ,
\]
where $V'(u)=\alpha u^{\alpha-1}$. It is obvious that $X$ is divergence-free with respect to the standard volume form $\Omega=dx_1\wedge\dots\wedge dx_n\wedge dp_1\wedge\dots\wedge dp_n$, that is $\Div_\Omega X=0$. A straightforward computation shows that the vector f\/ield
\[
Y_{ij}:=(x_ip_j-x_jp_i)\big[\alpha^{-1}(x_1\partial_{x_1}+\cdots+x_n\partial_{x_n})+p_1\partial_{p_1}+\dots+p_n\partial_{p_n}\big]
\]
verif\/ies
\[
[X,Y_{ij}]=\big(\alpha^{-1}-1\big)(x_ip_j-x_jp_i) X
\]
for all $1\leq i<j\leq n$. Applying Theorem~\ref{t1} with $f=1$ and $\lambda_{ij}=(\alpha^{-1}-1)(x_ip_j-x_jp_i)$ we conclude that $X$ has the following set of f\/irst integrals:
\[
H_{ij}=\Div_{\Omega} (Y_{ij})-\lambda_{ij}=\big(n\big(1+\alpha^{-1}\big)+2\big)(x_ip_j-x_jp_i) ,
\]
for all $1\leq i<j\leq n$, which just expresses the well known conservation law of the angular momentum. Observe that the f\/irst integral $H_{ij}$ is trivial only when $\alpha=\frac{-n}{n+2}$. Analogously, if we consider the vector f\/ield
\begin{gather*}
Y:=\left(\frac12\big(p_1^2+\dots+p_n^2\big)+V\big(x_1^2+\dots+x_n^2\big)\right)\\
\hphantom{Y:=}{}\times \big[\alpha^{-1}(x_1\partial_{x_1}+\dots+x_n\partial_{x_n})+
p_1\partial_{p_1}+\dots+p_n\partial_{p_n}\big]
\end{gather*}
we get that
\[
[X,Y]=\big(\alpha^{-1}-1\big)\left(\frac12\big(p_1^2+\dots+p_n^2\big)+V\big(x_1^2+\dots+x_n^2\big)\right) X ,
\]
and hence Theorem~\ref{t1} implies that $X$ has the following f\/irst integral as well:
\[
H=\big(n\big(1+\alpha^{-1}\big)+3-\alpha^{-1}\big)\left(\frac12\big(p_1^2+\dots+p_n^2\big)+V\big(x_1^2+\dots+x_n^2\big)\right) .
\]
This f\/irst integral ref\/lects the conservation of the energy, and is non-trivial unless $\alpha=\frac{1-n}{n+3}$. Note that this example includes the $n$-dimensional harmonic oscillator (with $\alpha=1$) and the $n$-dimensional Kepler problem (with $\alpha=(2-n)/2$) for $n>2$.
\end{example}

We f\/inish this letter by elaborating on the $2$-dimensional case.
This is a rather exceptional situation because the existence of a
normalizer implies the existence of an integrating factor, and vice
versa, and
therefore Theorem~\ref{t1} can be applied. In Example~\ref{Ex.2d}
below we show that the f\/irst integral $H$ in \eqref{e0} turns out to
be trivial for the usual choice of integrating factor (normalizer)
constructed from a given normalizer (an integrating factor).
Nevertheless, the f\/irst integral $H$ provided by Theorem~\ref{t1} is not trivial in general for $2$-dimensional vector f\/ields, as shown in Example~\ref{c1}.

\begin{example}\label{Ex.2d}
Let $X$ be a vector f\/ield on a $2$-dimensional manifold $M$ which
admits a transverse normalizer $Y$, i.e.\ $[X,Y]=\lambda X$ and
$\text{rank}(X,Y)=2$ a.e.\ in~$M$. It is immediate to
check~\cite{GG10} that the function $f:=(i_Xi_Y\Omega)^{-1}$ is an
integrating factor of $X$. We now compute the function $\lambda$ in
terms of $Y$ and $f$. Def\/ining the vector f\/ield $W:=fX$ and the
$1$-form $\omega:=i_W\Omega$, it follows that $i_Wi_Y\Omega=1$ and
using Cartan's formula for the Lie derivative we get
\begin{gather}\label{eq.ex1}
L_W(i_Y\Omega)=i_W(di_Y\Omega)+d(i_Wi_Y\Omega)=\Div_\Omega(Y) \omega .
\end{gather}
On the other hand, since the identity $L_Wi_Y=i_{[W,Y]}+i_YL_W$
holds, we derive that
\begin{gather}\label{eq.ex2}
L_W(i_Y\Omega)=\left(\lambda-\frac{L_Y(f)}{f}\right) \omega ,
\end{gather}
where we have used that $[W,Y]=\left(\lambda
-\frac{L_Y(f)}{f}\right) W$ and that $L_W(\Omega)=0$. Identifying
equations~\eqref{eq.ex1} and~\eqref{eq.ex2} and noticing that the
$1$-form $\omega$ is not trivial, we conclude that
\begin{gather*}
\lambda=\Div_\Omega(Y)+\frac{L_Y(f)}{f} .
\end{gather*}
Therefore, the f\/irst integral $H$ def\/ined in \eqref{e0} is
identically zero for this choice of integrating factor $f$ given a
normalizer $Y$.

Conversely, if $X$ admits an integrating factor $f$, it is easy to
check~\cite{GG10} that the vector f\/ield
$Y:=(f^2\Div_\Omega(X))^{-1} Z$ is a normalizer of $X$, with $Z$
the unique vector f\/ield which solves the equation $i_Z\Omega=-df$.
Arguing exactly as before it is straightforward to check the
normalizer condition $[X,Y]=\lambda X$, where
\[
\lambda=\Div_\Omega(Y)+\frac{L_Y(f)}{f} ,
\]
and therefore the f\/irst integral $H$ of \eqref{e0} is also trivial
for this choice of normalizer $Y$ given an integrating factor~$f$.

We observe that a non-trivial f\/irst integral of $X$ can be obtained through the classical Lie's integration algorithm~\cite{Ol93} for two dimensional vector f\/ields with a transverse normalizer $Y$. Indeed, the $1$-form
\[
\omega:=\frac{i_X\Omega}{i_Xi_Y\Omega}
\]
is closed, and hence there exists a smooth function $I$ def\/ined in any contractible domain of the manifold $M$ such that $dI=\omega$. It is clear that this function $I$ is a f\/irst integral of the vector f\/ield~$X$.
\end{example}

\begin{example}\label{c1}
Consider the following polynomial vector f\/ield in $\RR^2$ (we endow
the plane with Cartesian coordinates $(x,y)$ and the standard volume
form $dx\wedge dy$):
\[
X= \frac{1}{1+x^2}\left[-4y^3\partial_x+\big(1+x^4+y^4+4x^3\big)\partial_y\right] .
\]
It is easy to check that
\[
[X,Y]= \lambda X\,,
\]
where the vector f\/ield $Y$ is given by
\[
Y= \frac{(1+x^4+y^4)^2e^x}{\big(1+x^4+y^4+4x^3\big)^2+16y^6}\left[
\big(1+x^4+y^4+4x^3\big)\partial_x+4y^3\partial_y\right] ,
\]
and the expression of the function $\lambda$ is omitted for the sake
of simplicity. Moreover, the function
\[
f= (1+x^2)e^x
\]
is an integrating factor of the vector f\/ield $X$.

Therefore, we can apply Corollary~\ref{c0} to conclude that the
function
\[
H=\frac{L_Y(f)}{f}+\Div(Y)-\lambda= 2\big(1+x^4+y^4\big)e^x
\]
is a non-trivial global f\/irst integral of $X$.

For the sake of completeness, we f\/inally compute the f\/irst integral $I(x,y)$ of $X$ obtained through the classical Lie's integration algorithm (see the last paragraph of Example~\ref{Ex.2d}). After a few computations we get that
\[
I(x,y)=\int\frac{4y^3e^{-x}}{(1+x^4+y^4)^2}\,dy=-\frac{e^{-x}}{1+x^4+y^4}=-2H^{-1} ,
\]
thus showing that, in general, our approach gives a f\/irst integral dif\/ferent from the one obtained through classical methods, although of course in the two-dimensional case both f\/irst integrals are functionally dependent.
\end{example}

\subsection*{Acknowledgements}

The authors acknowledge three referees for corrections and comments which have helped to improve the presentation of this paper. J.Ll. is partially supported by a MICINN/FEDER grant no.
MTM2008--03437, by an AGAUR grant no. 2009SGR-410 and by ICREA
Academia. D.P.-S. is supported by a MICINN grant no. MTM2010--21186-C02-01, by the ICMAT Severo Ochoa grant no. SEV-2011-0087 and by the Ram\'on y Cajal program.

\pdfbookmark[1]{References}{ref}
\LastPageEnding

\end{document}